\UseRawInputEncoding
\documentclass[12pt]{article}
\usepackage{amssymb}
\usepackage{amsthm}
\usepackage{amsthm,amsmath,indentfirst,amssymb}
\usepackage{pstricks}
\usepackage{algorithm} %format of the algorithm
\usepackage{algorithmic} %format of the algorithm
\usepackage{blindtext}
\usepackage{graphicx}
\usepackage{cite}
\usepackage{float}
\usepackage[scriptsize,nooneline]{subfigure}
\usepackage{epstopdf}
\usepackage[colorlinks,
linkcolor=blue,
anchorcolor=blue,
citecolor=blue
]{hyperref}

{}
{}
{}

\begin{document}

\title{Necessary and sufficient conditions for boundedness of commutators of fractional integral operators on mixed-norm Lebesgue spaces
\thanks{The research was supported by the National Natural Science Fundation of China(12061069).}
}
\author{Houkun Zhang\thanks{E-mail address: zhanghkmath@163.com},\quad Jiang Zhou\thanks{ Corresponding author E-mail address: zhoujiang@xju.edu.cn.}
\\
\small  College of Mathematics and System Sciences, Xinjiang University, Urumqi 830046\\
\small People's Republic of China
}
\date{}
\maketitle
{\bf Abstract:} In this paper, the sharp maximal theorem is generalized to mixed-norm ball Banach function spaces, which is defined as Definition 2.7. As an application, we give a characterization of BMO via the boundedness of commutators of fractional integral operators on mixed-norm Lebesgue spaces. Moreover, the characterization of homogeneous Lipschitz space is also given by the boundedness of commutators of fractional integral operators on mixed-norm Lebesgue spaces. Finally, two applications of Corollary 6.4 are given.

{\bf MSC(2000) subject classification:}  42B25; 42B20.

\maketitle

\section{Introduction}\label{sec1}
\par

In recent years, due to the more precise structure of mixed-norm function spaces than the corresponding classical function spaces, the mixed-norm function spaces are widely used in the partial differential equations \cite{3,4,5,6}. In 1961, the mixed Lebesgue spaces $L^{\vec{p}}(\mathbb{R}^n)(0<\vec{p}<\infty)$ were studied by Benedek and Panzone \cite{7}. These spaces were
natural generalizations of the classical Lebesgue spaces $L^p(0<p<\infty)$. After that, many function spaces with mixed norm were introduced, such as mixed-norm Lorentz spaces \cite{8}, mixed-norm Lorentz-Marcinkiewicz spaces \cite{9}, mixed-norm Orlicz spaces \cite{10}, anisotropic mixed-norm Hardy spaces \cite{11}, mixed-norm Triebel-Lizorkin spaces \cite{12}, mixed Morrey spaces \cite{15,16} and weak mixed-norm Lebesgue spaces \cite{13}. More information can be found in \cite{14}.

In 2017, ball (quasi-)Banach function spaces were introduced by Sawano et al. \cite{19}. Ball (quasi-)Banach function spaces are the generalizations of the (quasi-)Banach function spaces (see Remark 2.6).
In this paper, according to the definitions of the mixed-norm Banach function spaces \cite{44} defined by Blozinski, we introduce the definitions of mixed-norm ball (quasi-)Banach function spaces and prove sharp maximal theorem mixed-norm ball Banach function spaces (see Theorem 3.7). we also acquire Corollary 3.8 and corollary 3.9.

Moreover, we point out that many spaces are ball (quasi-)Banach function spaces, such as variable Lebesgue spaces \cite{32} and Orlicz spaces \cite{33}, which also are (quasi-)Banach function spaces,  Besides, some spaces are ball (quasi-)Banach function spaces, which are not necessary to be (quasi-)Banach function spaces, such as Morrey spaces \cite{31}, mixed-norm Lebesgue spaces, weighted Lebesgue spaces, and Orlicz-slice spaces \cite{35}.

Particularly, mixed-norm Lebesgue spaces also can be regarded as mixed-norm ball (quasi-)Banach function spaces. In addition, mixed-norm Lorentz spaces, mixed-norm Orlicz spaces, and mixed-norm Lebesgue spaces with variable exponents are also mixed-norm ball (quasi-)Banach function spaces. For more studies of ball Banach function spaces, we refer the readers to \cite{37,38,39,40,41}.

Given a number $0<\alpha<n$, the fractional integral functions of a measurable function $f$ on $\mathbb{R}^n$ is defined by
$$I_\alpha f(x)=\int_{\mathbb{R}^n}\frac{f(y)}{|x-y|^{n-\alpha}}dy.$$
These operators $I_\alpha$ play an essential role in real and harmonic analysis \cite{1,2}.

For a locally integrable function $b$, the commutator is defined by
$$[b,I_\alpha]f(x):=b(x)I_\alpha f(x)-I_\alpha(bf)(x)=\int_{\mathbb{R}^n}\frac{(b(x)-b(y))f(y)}{|x-y|^{n-\alpha}}dy,$$
which are introduced by Chanillo in \cite{17}. It is obvious that
$$|[b,I_\alpha](f)(x)|\le I_{\alpha,b}(|f|)(x),$$
where
$$I_{\alpha,b}(f)(x):=I_{\alpha}(|b(x)-b(\cdot)|f)(x)=\int_{\mathbb{R}^n}\frac{|b(x)-b(y)|f(y)}{|x-y|^{n-\alpha}}dy.$$

Many classical works on the characterizations of $BMO(\mathbb{R}^n)$ and homogeneous Lipschitz spaces were studied via the boundedness of $[b,I_\alpha]$ on classical Lebesgue spaces. In 1982, an early characterization of $BMO(\mathbb{R}^n)$ spaces was investigated by Chanillo \cite{42} via the $(L^p(\mathbb{R}^n),L^q(\mathbb{R}^n))$-boundedness of $[b,I_\alpha]$. That is
$$b\in BMO(\mathbb{R}^n)\Leftrightarrow [b,I_\alpha]:L^p(\mathbb{R}^n)\mapsto L^q(\mathbb{R}^n),$$
where $1<p<q<\infty,~\frac{1}{p}-\frac{1}{q}=\frac{\alpha}{n}$. the characterization of homogeneous Lipschitz spaces $Lip_{\beta}(\mathbb{R}^n)(0<\beta<1)$ was given by Paluszynski \cite{43} via $(L^p(\mathbb{R}^n),L^q(\mathbb{R}^n))$-boundedness of $[b,I_\alpha]$ as:
$$b\in Lip_{\beta}(\mathbb{R}^n)\Leftrightarrow [b,I_\alpha]:L^p(\mathbb{R}^n)\mapsto L^q(\mathbb{R}^n),$$
where $1<p<q<\infty,~\frac{1}{p}-\frac{1}{q}=\frac{\alpha+\beta}{n}$.

In 2019, Nogayama given a characterization of $BMO(\mathbb{R}^n)$ via the boundedness of commutators of fractional integral operators on mixed Morrey spaces \cite{16}. Due to the definitions of mixed Morrey spaces, we know that
$$b\in BMO(\mathbb{R}^n)\Leftrightarrow [b,I_\alpha]:L^{\vec{p}}(\mathbb{R}^n)\mapsto L^{\vec{q}}(\mathbb{R}^n),$$
where $1<\vec{p},\vec{q}<\infty$ and
$$\alpha=\sum_{i=1}^n\frac{1}{p_i}-\sum_{i=1}^n\frac{1}{q_i},~p_j\sum_{i=1}^n\frac{1}{p_i}=q_j\sum_{i=1}^n\frac{1}{q_i}~~(j=1,\cdots,n).$$

In fact, the condition
$$p_j\sum_{i=1}^n\frac{1}{p_i}=q_j\sum_{i=1}^n\frac{1}{q_i}~~(j=1,\cdots,n)$$
is not necessary. In Theorem 4.3, a weaker condition is given. Furthermore, the characterization of homogeneous Lipschitz spaces $Lip_{\beta}(\mathbb{R}^n)(0<\beta<1)$ was also given in Theorem 5.3 via the boundedness of $[b,I_\alpha]$ on mixed-norm Lebesgue spaces.

The paper is organized as the following. In Section 2, some notations and notions are recalled and we introduce the definitions of mixed-norm ball (quasi-)Banach function spaces. In Section 3, the sharp maximal theorem is investigated on mixed-norm ball Banach function spaces via corresponding extrapolation theorem. As an application, we give a characterization of BMO via the boundedness of commutators of fractional integral operators on Mixed-norm Lebesgue spaces in Section 4. Moreover, the characterization of homogeneous Lipschitz spaces is also given in Section 5. In Section 6, two applications of Corollary 6.4 are given.

\section{Some notations and notions}\label{sec2}

In this section, we make some conventions on notation and recall some notions. Let $\vec{p}=(p_1,p_2,\cdots,p_n),~\vec{q}=(q_1,q_2,\cdots,q_n)$, are n-tuples and $1<p_i,q_i<\infty,~i=1,2,\cdots,n$. We define that  if $\varphi(a,b)$ is a relation or equation among numbers, $\varphi(\vec{p},\vec{q})$ will mean that $\varphi(p_i,q_i)$ hords for each $i$. For example, $\vec{p}<\vec{q}$ means that $p_i<q_i$ holds for each $i$ and $\frac{1}{\vec{p}}+\frac{1}{\vec{p}\,'}=1$ means $\frac{1}{p_i}+\frac{1}{p_i'}=1$ hold for each $i$. The symbol $Q$ denote the cubes whose edges are parallel to the coordinate axes and $Q(x,r)$ denote the open cube centered at $x$ of side length $r$. Let $cQ(x,r)=Q(x,cr)$. Denote by the symbol $\mathfrak{M}(\mathbb{R}^n)$ the set of all measurable function on $\mathbb{R}^n$. $A\thickapprox B$ means that $A$ is equivalent to $B$. That is $A\le CB$ and $B\le CA$, where $C$ is a constant. Through all paper, every constant $C$ is not necessarily equal. Now, let us recall some definitions.

The definitions of some maximal functions and some class of weight functions are given as the following.

\textbf{Definition 2.1} Given a locally integrable function $f$. The Hardy-Littlewood maximal operators is defined by
$$M(f)(x):=\sup_{Q\ni x}\frac{1}{|Q|}\int_{Q}|f(y)|dy,$$
where the supremum is taken over all cube $Q\subset\mathbb{R}^n$ containing $x$. The sharp maximal operators is defined by
$$M^{\sharp}(f)(x):=\sup_{Q\ni x}\frac{1}{|Q|}\int_{Q}|f(y)-f_Q|dy,$$
where $f_Q=\frac{1}{|Q|}\int_{Q}f$ and the supremum is taken over all cube $Q\subset\mathbb{R}^n$ containing $x$.

\textbf{Definition 2.2} An $A_p(\mathbb{R}^n)$ weight $\omega$, with $1\le p<\infty$, is a locally integrable and nonnegative function on $\mathbb{R}^n$ satisfying that, when $1<p<\infty$,
$$[\omega]_{A_p(\mathbb{R}^n)}:=\sup_{Q\subset\mathbb{R}^n}\left(\frac{1}{|Q|}\int_Q\omega(x)dx\right) \left(\frac{1}{|Q|}\int_Q\omega(x)^{\frac{1}{1-p}}dx\right)^{p-1}<\infty,$$
and, when $p=1$,
$$[\omega]_{A_1(\mathbb{R}^n)}:=\sup_{Q\subset\mathbb{R}^n}\left(\frac{1}{|Q|}\int_Q\omega(x)dx \cdot\|\omega^{-1}\|_{L^{\infty}(\mathbb{R}^n)}\right)<\infty.$$
Define $A_{\infty}(\mathbb{R}^n):=\bigcup_{1\le p<\infty}A_p(\mathbb{R}^n)$. It is well-know that $A_p(\mathbb{R}^n)\subset A_q(\mathbb{R}^n)$ for $1\le p\le q\le\infty$.

The following is the corresponding weight functions for the product domain $\mathbb{R}^{n_1}\times\mathbb{R}^{n_2}$ from Chapter $\text{IV}$ of \cite{29}.

\textbf{Definition 2.3} Let $Q_1\subset\mathbb{R}^{n_1}$, $Q_2\subset\mathbb{R}^{n_2}$, $R=Q_1\times Q_2$ and $n=n_1+n_2$. For $1<p<\infty$, a locally integrable and nonnegative function on $\mathbb{R}^n$ is said to be an $A_p^{*}(\mathbb{R}^n)$ weight if
$$[\omega]_{A_p^{*}(\mathbb{R}^n)}:=\sup_{R\subset\mathbb{R}^n}\left(\frac{1}{|R|}\int_R\omega(z)dz\right) \left(\frac{1}{|R|}\int_R\omega(z)^{\frac{1}{1-p}}dz\right)^{p-1}<\infty.$$
A locally integrable and nonnegative function on $\mathbb{R}^n$ is said to be an $A_1^{*}(\mathbb{R}^n)$ weight if
$$[\omega]_{A^{*}_1(\mathbb{R}^n)}:=\sup_{R\subset\mathbb{R}^n}\left(\frac{1}{|R|}\int_R\omega(z)dz \cdot\|\omega^{-1}\|_{L^{\infty}(\mathbb{R}^n)}\right)<\infty.$$
We write $A^{*}_{\infty}(\mathbb{R}^n):=\bigcup_{1\le p<\infty}A^{*}_p(\mathbb{R}^n)$.

By the definition of $A_p^{*}(\mathbb{R}^n)$ for $1\le p<\infty$, it is easy to prove that if $\mu\in A_p(\mathbb{R}^{n_1})$ and $\nu\in A_p(\mathbb{R}^{n_2})$ then
$$\mu\nu\in A_p^{*}(\mathbb{R}^n)$$
where $n=n_1+n_2$. Moreover, $[\omega]_{A_p(\mathbb{R}^n)}\le [\omega]_{A_p^{*}(\mathbb{R}^n)}$ and $A_p^{*}(\mathbb{R}^n)\subset A_p(\mathbb{R}^n)$ for $1\le p<\infty$.

\textbf{Definition 2.4} Let $0<p<\infty$ and $\omega\in A_{\infty}(\mathbb{R}^n)$. The weighted Lebesgue space $L_{\omega}^p(\mathbb{R}^n)$ is defined to be the set of all measurable functions $f$ on $\mathbb{R}^n$ such that
$$\|f\|_{L_{\omega}^p(\mathbb{R}^n)}:=\left(\int_{\mathbb{R}^n}|f(x)|^p\omega(x)dx\right)^{\frac{1}{p}}<\infty.$$

The definition of ball (quasi-)Banach function spaces is presented as the following, which were introduced by Sawano et al. \cite{19}.

\textbf{Definition 2.5} A (quasi-)Banach space $\mathbf{X}\subset\mathfrak{M}(\mathbb{R}^n)$ with (quasi-)norm $\|\cdot\|_{\mathbf{X}}$ is called a ball (quasi-)Banach function space if it satisfies\\
(1) $|g|\le|f|$ almost everywhere implies that $\|g\|_{\mathbf{X}}\le\|f\|_{\mathbf{X}}$;\\
(2) $0\le f_m\uparrow f$ almost everywhere implies that $\|f_m\|_{\mathbf{X}}\uparrow\|f\|_{\mathbf{X}};$\\
(3) If $|Q|<\infty$, then $\chi_Q\in \mathbf{X}$;\\
(4) If $f\ge 0$ almost everywhere and $|Q|<\infty$, then
$$\int_Q f(x)dx\le C_Q\|f\|_{\mathbf{X}};$$
for some constant $C_Q$, $0<c_Q<\infty$, depending on $Q$ but independent of $f$.

\textbf{Remark 2.6} The definition remains unchanged if we replace "cube" with "ball" in the above. So this definition deserves this name. Particularly, if we replace any cube $Q$ by any measurable set $E$ in Definition 2.5 and add the condition $\|f\|_{\mathbf{X}}=\||f|\|_{\mathbf{X}}$, we obtain the definition of (quasi-)Banach function space (see Definition 1.1 of Chapter 1 of \cite{36}).

Next, we given the definition of mixed-norm ball (quasi-)Banach function spaces.

\textbf{Definition 2.7} Let $n_1,~n_2\in \mathbb{N}$ and $\mathbf{X}_1$ and $\mathbf{X}_2$ be ball (quasi-)Banach function spaces on $\mathbb{R}^{n_1}$ and $\mathbb{R}^{n_2}$, respectively. The mixed-norm ball (quasi-)Banach function spaces $(\mathbf{X}_1,\mathbf{X}_2)$ consists of Lebesgue measurable function on $\mathbb{R}^{n_1}\times\mathbb{R}^{n_2}$, $f$ such that
$$\|f\|_{(\mathbf{X}_1,\mathbf{X}_2)}=\|\|f\|_{\mathbf{X}_1}\|_{\mathbf{X}_2}<\infty.$$

It is easy to prove that $(\mathbf{X}_1,\mathbf{X}_2)$ is also a ball (quasi-)Banach function space. In fact, if $Q=Q_1\times Q_2\in\mathbb{R}^{n_1}\times\mathbb{R}^{n_2}$, then
$$\|\|\chi_Q\|_{\mathbf{X}_1}\|_{\mathbf{X}_2}=\|\chi_{Q_1}\|_{\mathbf{X}_1}\|\chi_{Q_2}\|_{\mathbf{X}_2}<\infty$$
and
\begin{align*}
\int_{Q}f(x,y)dxdy&=\int_{Q_2}\int_{Q_1}f(x,y)dxdy\\
&\le C_{Q_1}\int_{Q_2}\|f(\cdot,y)\|_{\mathbf{X}_1}dy\\
&\le C_{Q}\|\|f\|_{\mathbf{X}_1}\|_{\mathbf{X}_2}.
\end{align*}

The definition of associate space of a ball Banach function space can be found in chapter 1 of \cite{22}. The definition of associate space of a ball (quasi-)Banach function space is given as the following.

\textbf{Definition 2.8} For any ball (quasi-)Banach function space $\mathbf{X}$, the associate space (also called the K\"othe dual) $\mathbf{X}'$ is defined by setting
$$\mathbf{X}':=\left\{f\in \mathfrak{M}(\mathbb{R}^n):\|f\|_{\mathbf{X}'}: =\sup_{g\in\mathbf{X},\|g\|_{\mathbf{X}}=1}\int_{\mathbb{R}^n}|f(x)g(x)|dx<\infty\right\},$$
where $\|\cdot\|_{\mathbf{X}'}$ is called the associate norm of $\|\cdot\|_{\mathbf{X}}.$

We still need to recall the notion of the convexity of ball (quasi-)Banach function spaces, which can be found in Definition 4.6 of \cite{19}.

\textbf{Definition 2.9} Let $\mathbf{X}$ be a ball (quasi-)Banach function space and $0<p<\infty$. The p-convexification $\mathbf{X}^p$ of $\mathbf{X}$ is defined by setting
$$\mathbf{X}^p:=\left\{f\in\mathfrak{M}(\mathbb{R}^n):|f|^p\in\mathbf{X}\right\}$$
equipped with the (quasi-)norm $\|f\|_{\mathbf{X}^p}:=\| |f|^p\|_{\mathbf{X}}^{\frac{1}{p}}$.

\textbf{Remark 2.10} (1) By the definition 2.9, we know that
$$\|f\|_{(\mathbf{X}_1,\mathbf{X}_2)^p}=\|\||f|^p\|_{\mathbf{X}_1}\|_{\mathbf{X}_2}^{\frac{1}{p}} =\|\|f\|_{\mathbf{X}_1^p}\|_{\mathbf{X}_2^p}.$$
Thus,
$$(\mathbf{X}_1,\mathbf{X}_2)^p=(\mathbf{X}_1^p,\mathbf{X}_2^p).\eqno{(2.1)}$$

(2) If $f\in (\mathbf{X}_1,\mathbf{X}_2)'$, then there exist $h_1\in\mathbf{X}_1$, $h_2\in\mathbf{X}_2$, $\|h_1\|_{\mathbf{X}_1}=1$ and $\|h_2\|_{\mathbf{X}_2}=1$ such that
\begin{align*}
\|f\|_{(\mathbf{X}'_1,\mathbf{X}'_2)}&=\|\|f\|_{\mathbf{X}'_1}\|_{\mathbf{X}'_2}\\
&\le 2\int_{\mathbb{R}^{n_1}\times\mathbb{R}^{n_2}}|f(x,y)h_1(x)h_2(y)|dxdy\\
&\le \sup_{g\in(\mathbf{X}_1,\mathbf{X}_2),\|g\|_{(\mathbf{X}_1,\mathbf{X}_2)}=1} 2\int_{\mathbb{R}^{n_1}\times\mathbb{R}^{n_2}}|f(x,y)g(x,y)|dxdy\\
&\le \|f\|_{(\mathbf{X}_1,\mathbf{X}_2)'}
\end{align*}
For another hand, if $f\in (\mathbf{X}'_1,\mathbf{X}'_2)$, then there exist $g\in(\mathbf{X}_1,\mathbf{X}_2)$,
\begin{align*}
\|f\|_{(\mathbf{X}_1,\mathbf{X}_2)}&\le 2\int_{\mathbb{R}^{n_1}\times\mathbb{R}^{n_2}}|f(x,y)g(x,y)|dxdy\\
&\le 2\int_{\mathbb{R}^{n_2}}\|f(\cdot,y)\|_{\mathbf{X}'_1}\|g(\cdot,y)\|_{\mathbf{X}_1}dxdy\\
&\le 2\|f\|_{(\mathbf{X}'_1,\mathbf{X}'_2)}\|g\|_{(\mathbf{X}_1,\mathbf{X}_2)}\\
&= 2\|f\|_{(\mathbf{X}_1,\mathbf{X}_2)'}
\end{align*}
Hence, $(\mathbf{X}'_1,\mathbf{X}'_2)$ are equivalent to $(\mathbf{X}_1,\mathbf{X}_2)'$.

Obviously, if $\mathbf{X}$ is a ball (quasi-)Banach function space, the $\mathbf{X}^p$ and $\mathbf{X}'$ are also ball (quasi-)Banach function spaces. Now, the definitions of mixed-norm Lebesgue spaces are given as the following, which were introduced by Benedek and Panzone \cite{7}.

\textbf{Definition 2.11} Let $f$ is a measurable function on $\mathbb{R}^n$ and $1<\vec{p}<\infty$. We say that $f$ belongs to the mixed Lebesgue spaces $L^{\vec{p}}(\mathbb{R}^n)$, if the norm
$$\left\|f\right\|_{L^{\vec{p}}(\mathbb{R}^n)}=\left(\int_{\mathbb{R}}\cdots\left(\int_{\mathbb{R}}\left|f(x)\right|^{p_1}\,dx_1\right)
^{\frac{p_2}{p_1}}\cdots\,dx_n\right)^{\frac{1}{p_n}}<\infty.$$
Note that if $p_1=p_2=\cdots=p_n=p$, then $L^{\vec{p}}(\mathbb{R}^n)$ are reduced to classical Lebesgue spaces $L^p$ and
$$\left\|f\right\|_{L^{\vec{p}}(\mathbb{R}^n)}=\left(\int_{\mathbb{R}^n}\left|f(x)\right|^{p} dx\right)^{\frac{1}{p}}.$$

\textbf{Remark 2.12} It is easy to calculate that for any cube $Q\subset\mathbb{R}^n$ and $|Q|<\infty$,
$$\|\chi_{Q}\|_{L^{\vec{p}}(\mathbb{R}^n)}=|Q|^{\frac{1}{n}\sum_{i=1}^{n}\frac{1}{p_i}}$$
and
$$\int_{Q}f(x)dx\le |Q|^{\frac{1}{n}\sum_{i=1}^{n}\frac{1}{p'_i}}\|f\|_{L^{\vec{p}}(\mathbb{R}^n)}.$$
By Levi Lemma, (2) also can be proved. Thus $L^{\vec{p}}$ are ball Banach function spaces. But it is uncertain wether mixed-norm Lebesgue spaces are Banach function spaces.

Let us recall the definition of $BMO(\mathbb{R}^n)$ spaces.

\textbf{Definition 2.13} If $b$ is a measurable function on $\mathbb{R}^n$ and satisfies that
$$\|b\|_{BMO(\mathbb{R}^n)}=\sup_{Q\subset\mathbb{R}^n}\frac{1}{|Q|}\int_{Q}|b(y)-b_Q|dy<\infty,$$
then $b\in BMO(\mathbb{R}^n)$ and $\|b\|_{BMO(\mathbb{R}^n)}$ are the norms of $b$ in $BMO(\mathbb{R}^n)$.

\section{Sharp maximal theorem on mixed-norm ball quasi-Banach function spaces}\label{sec3}
\par

In this section, we will prove an extrapolation theorem on mixed-norm ball quasi-Banach function spaces. According to the extrapolation theorem, we prove the sharp maximal theorem on mixed-norm ball quasi-Banach function spaces. The sharp maximal theorem on weighted Lebesgue spaces can be found in Theorem 3.4.5 of \cite{34}.

\textbf{Lemma 3.1} Let $1<p<\infty$, $\omega\in A_p$. Then for any $f\in L_\omega^p$,
$$\int_{\mathbb{R}^n}\left(Mf(x)\right)^p\omega(x)dx\le C\int_{\mathbb{R}^n}\left(M^{\sharp} f(x)\right)^p\omega(x)dx$$
holds.

To prove the sharp maximal theorem on ball (quasi-)Banach function space, we give an assumption and some lemmas.

\textbf{Assumption 3.2} There exists an $s\in(1,\infty)$ such that $\mathbf{X}$ is a ball Banach function space, and that $Mf$ is bounded on $(\mathbf{X}^{\frac{1}{s}})'$.

The following result can be found in Lemma 4.7 of \cite{21}.

\textbf{Lemma 3.3} Let $\mathbf{X}$ be a ball quasi-Banach function space satisfying Assumption 2.10. Then there exists an $0<\varepsilon<1$ such that $\mathbf{X}$ continuously embeds into $L^s_{\omega}(\mathbb{R}^n)$ with $\omega:=[M(\chi_{Q(0,1)})]^{\varepsilon}\in A_1(\mathbb{R}^n)$, namely, there exists a positive constant $C$ such that, for any $f\in\mathbf{X}$,
$$\|f\|_{L_{\omega}^s(\mathbb{R}^n)}\le C\|f\|_{\mathbf{X}}.$$

Due to Lemma 3.3, the following lemma is proved.

\textbf{Lemma 3.4} Let $\mathbf{X}_1$ and $\mathbf{X}_2$ be ball quasi-Banach function spaces on $\mathbb{R}^{n_1}$ and $\mathbb{R}^{n_2}$, and satisfy Assumption 2.10. Then there exists an $0<\varepsilon<1$ such that $(\mathbf{X}_1,\mathbf{X}_2)$ continuously embeds into $L^s_{\omega}(\mathbb{R}^{n_1+n_2})$ with $\omega=\omega_1\omega_2\in A_1^{*}(\mathbb{R}^n)$, where $\omega_1:=[M(\chi_{Q(0,1)})]^{\varepsilon}\in A_1(\mathbb{R}^{n_1})$ $\omega_2:=[M(\chi_{Q(0,1)})]^{\varepsilon}\in A_1(\mathbb{R}^{n_2})$ and $n=n_1+n_2$, namely, there exists a positive constant $C$ such that, for any $f\in(\mathbf{X}_1,\mathbf{X}_2)$,
$$\|f\|_{L_{\omega}^s(\mathbb{R}^{n_1+n_2})}\le C\|f\|_{(\mathbf{X}_1,\mathbf{X}_2)}.$$

\textbf{Proof} If $f\in(\mathbf{X}_1,\mathbf{X}_2)$, then by Lemma 3.3
\begin{align*}
\|f\|_{L_{\omega}^s(\mathbb{R}^{n_1+n_2})}&=\left(\int_{\mathbb{R}^{n_2}} \left(\int_{\mathbb{R}^{n_1}}|f(x,y)|^s\omega_1(x)dx\right)^{\frac{1}{s}\cdot s}\omega_2(y)dy\right)^{\frac{1}{s}}\\
&\le\left(\int_{\mathbb{R}^{n_2}}\|f(\cdot,y)\|^s_{\mathbf{X}_1}\omega_2(y)dy\right)^{\frac{1}{s}}\\
&\le\|\|f\|_{\mathbf{X}_1}\|_{\mathbf{X}_2}
\end{align*}
The proof is completed. $~~~~\blacksquare$

The following extrapolation theorem plays an important role in the proof of Theorem 3.4. The extrapolation theorem is a slight variant of a special case of Theorem 3.2 of \cite{30}, via replacing Banach function spaces by ball quasi-Banach function spaces. The proof of Lemma 3.5 is similar to Theorem 3.2 of \cite{30}.

\textbf{Lemma 3.5} Let $\mathbf{X}_1$ and $\mathbf{X}_2$ be a ball quasi-Banach function spaces and $p_0\in(0,\infty)$. Let $\mathcal{F}$ be the set of all pairs of nonnegative measurable functions $(F,G)$ such that, for any given $\omega\in A^{*}_1(\mathbb{R}^n)$,
$$\int_{\mathbb{R}^{n_1}\times\mathbb{R}^{n_2}}\left(F(x,y)\right)^{p_0}\omega(x,y)dxdy\le C_{(p_0,[\omega]_{A^{*}_1(\mathbb{R}^n)})}\int_{\mathbb{R}^{n_1}\times\mathbb{R}^{n_2}}\left(G(x,y)\right)^{p_0}\omega(x,y)dxdy.$$
where $C_{(p_0,[\omega]_{A^{*}_1(\mathbb{R}^n)})}$ is a positive constant independent of $(F,G)$, but dependents on $p_0$ and $A^{*}_1(\mathbb{R}^n)$. Assume that there exisitive a $q_0\in[p_0,\infty)$ such that $\mathbf{X}_1^{\frac{1}{q_0}}$ and $\mathbf{X}_2^{\frac{1}{q_0}}$ are ball Banach function spaces and $Mf$ is bounded on $(\mathbf{X}_1^{\frac{1}{q_0}})'$ and $(\mathbf{X}_2^{\frac{1}{q_0}})'$. Then there exists a positive constant $C_0$ such that, for any $(F,G)\in\mathcal{F}$,
$$\|F\|_{(\mathbf{X}_1,\mathbf{X}_2)}\le C_0\|G\|_{(\mathbf{X}_1,\mathbf{X}_2)}.$$

For completeness, we will give the proof of Lemma 3.5. Before that, the following lemma is necessary. It can be found in Theorem 3.1 of \cite{30}.

\textbf{Lemma 3.6} Suppose that for $P_0$ with $1\le p_0<\infty$. Let $\mathcal{F}$ be the set of all pairs of nonnegative measurable functions $(F,G)$ such that, for any given $\omega\in A^{*}_1(\mathbb{R}^n)$,
$$\int_{\mathbb{R}^{n_1}\times\mathbb{R}^{n_2}}F(x,y)^{p_0}\omega(x,y)dxdy\le C\int_{\mathbb{R}^{n_1}\times\mathbb{R}^{n_2}}G(x,y)^{p_0}\omega(x,y)dxdy~~~~(F,G)\in\mathcal{F},$$
where $n=n_1+n_2$. Then, for all $1<p<\infty$ and $\omega\in A^{*}_p(\mathbb{R}^n)$, we have
$$\int_{\mathbb{R}^{n_1}\times\mathbb{R}^{n_2}}F(x,y)^{p}\omega(x,y)dxdy\le C\int_{\mathbb{R}^{n_1}\times\mathbb{R}^{n_2}}G(x,y)^{p}\omega(x,y)dxdy~~~~(F,G)\in\mathcal{F}.$$

\textbf{The proof of Lemma 3.5} According to Theorem 3.6, if let
$$\mathcal{F}_0=\{(F^{p_0},G^{p_0}):(F,G)\in\mathcal{F}\}$$
then we find that for any $p\ge 1$ and $\omega\in A^{*}_1(\mathbb{R}^n)\subset A^{*}_p(\mathbb{R}^n)$,
$$\int_{\mathbb{R}^{n_1}\times\mathbb{R}^{n_2}}F(x,y)^{p_0p}\omega(x,y)dxdy\le C\int_{\mathbb{R}^{n_1}\times\mathbb{R}^{n_2}}G(x,y)^{p_0p}\omega(x,y)dxdy.$$
Let $p=\frac{q_0}{q_0}$. Then
$$\int_{\mathbb{R}^{n_1}\times\mathbb{R}^{n_2}}F(x,y)^{q_0}\omega(x,y)dxdy\le C\int_{\mathbb{R}^{n_1}\times\mathbb{R}^{n_2}}G(x,y)^{q_0}\omega(x,y)dxdy.\eqno{(3.1)}$$

We write $Y_1=\mathbf{X}_1^{\frac{1}{q_0}}$ and $Y_2=\mathbf{X}_2^{\frac{1}{q_0}}$. For any functions $h_1(x)\in\mathbf{Y}'_1$ and $h_2(y)\in\mathbf{Y}'_2$, let
$$R_1h_1(x)=\sum_{k=0}^{\infty}\frac{M^kh_1(x)}{(2A_1)^k}$$
and
$$R_2h_2(y)=\sum_{k=0}^{\infty}\frac{M^kh_2(y)}{(2A_2)^k},$$
where $A_1=\|M\|_{\mathbf{Y}'_1\rightarrow\mathbf{Y}'_1}$, $A_2=\|M\|_{\mathbf{Y}'_2\rightarrow\mathbf{Y}'_2}$, $M^0h=|h|$ and $M^kh=M(M^{k-1}h)$. Obviously, $h_1$ and $h_2$ satisfy that
$$|h_i|\le R_ih_i,\eqno{(3.2)}$$
$$\|R_ih_i\|_{\mathbf{Y}'_i}\le 2\|h_i\|_{\mathbf{Y}'_i},\eqno{(3.3)}$$
$$M(R_ih_i)\le 2A_iR_ih_i,\eqno{(3.4)}$$
for $i=1,2$. (3.2), (3.3) and (3.4) follow from the definition of $R_1$ and $R_2$. By (3.4), it is obvious that
$$R_1h_1(x)\in A_1(\mathbb{R}^{n_1}),~R_2h_2(y)\in A_1(\mathbb{R}^{n_2})$$
and
$$R_1h_1(x)R_2h_2(y)\in A^{*}_1(\mathbb{R}^{n_1+n_2}).\eqno{(3.5)}$$

According to the definition of associate spaces and (2.1) there exists measurable functions $h_1(x)\in\mathbf{Y}'_1$ and $h_2(y)\in\mathbf{Y}'_2$ such that $\|h_1\|_{\mathbf{Y}'_1}=1$, $\|h_2\|_{\mathbf{Y}'_2}=1$ and
\begin{align*}
\|F\|^{q_0}_{(\mathbf{X}_1,\mathbf{X}_2)}&=\||F|^{q_0}\|_{(\mathbf{Y}_1,\mathbf{Y}_2)}\\
&\le C\int_{\mathbb{R}^{n_2}}\||F(\cdot,y)|^{q_0}\|_{\mathbf{Y}_1}|h_2(y)|dy\\
&\le C\int_{\mathbb{R}^{n_1}\times\mathbb{R}^{n_2}}|F(x,y)|^{q_0}|h_1(x)||h_2(y)|dxdy
\end{align*}
According to (3.2), we have
$$\|F\|^{q_0}_{(\mathbf{X}_1,\mathbf{X}_2)}\le C\int_{\mathbb{R}^{n_1}\times\mathbb{R}^{n_2}}|F(x,y)|^{q_0}|R_1h_1(x)||R_2h_2(y)|dxdy.$$
Due to (3.1) and (3.5), we have
$$\|F\|^{q_0}_{(\mathbf{X}_1,\mathbf{X}_2)}\le C\int_{\mathbb{R}^{n_1}\times\mathbb{R}^{n_2}}|G(x,y)|^{q_0}|R_1h_1(x)||R_2h_2(y)|dxdy.$$
By the definition of associate spaces, H\"older's inequalities are obtained. That is
$$\int_{\mathbb{R}^{n}}f(x)g(y)dxdy\le \|f\|_{\mathbf{X}}\|g\|_{\mathbf{X}'}$$
for $f\in\mathbf{X}$ and $g\in\mathbf{X}'$. Thus,
\begin{align*}
\|F\|^{q_0}_{(\mathbf{X}_1,\mathbf{X}_2)}&\le C\int_{\mathbb{R}^{n_2}}\||G(\cdot,y)|^{q_0}\|_{\mathbf{Y}_1}\|R_1h_1\|_{\mathbf{Y}'_1}|R_2h_2(y)|dy\\
&\le C\||G|^{q_0}\|_{(\mathbf{Y}_1,\mathbf{Y}_2)}\|R_1h_1\|_{\mathbf{Y}'_1}\|R_2h_2\|_{\mathbf{Y}'_2}
\end{align*}
By (3.3),
\begin{align*}
\|F\|^{q_0}_{(\mathbf{X}_1,\mathbf{X}_2)}
&\le C\||G|^{q_0}\|_{(\mathbf{Y}_1,\mathbf{Y}_2)}\|h_1\|_{\mathbf{Y}'_1}\|h_2\|_{\mathbf{Y}'_2}\\
&=C\||G|^{q_0}\|_{(\mathbf{Y}_1,\mathbf{Y}_2)}\\
&=C\|G\|^{q_0}_{(\mathbf{X}_1,\mathbf{X}_2)}.
\end{align*}
The proof is completed. $~~~~\blacksquare$

\textbf{Theorem 3.7}  Let $\mathbf{X}_1$ and $\mathbf{X}_2$ be a ball quasi-Banach function spaces and satisfy Assumption 3.2 for $1<s_0<\infty$. Then there exists a positive constant $C_0$ such that, for any $f\in(\mathbf{X}_1,\mathbf{X}_2)$,
$$\|Mf\|_{(\mathbf{X}_1,\mathbf{X}_2)}\le C\|M^{\sharp} f\|_{(\mathbf{X}_1,\mathbf{X}_2)}.$$

\textbf{Proof} Let
$$\mathcal{F}:=\left\{(Mf,M^{\sharp}f):f\in\bigcup_{\omega\in A^{*}_1(\mathbb{R}^n)}L_{\omega}^p(\mathbb{R}^n)\right\}.$$
According to Lemma 3.1,
$$\int_{\mathbb{R}^n}\left(Mf(x)\right)^p\omega(x)dx\le C_{(p,[\omega]_{A^{*}_1(\mathbb{R}^n})}\int_{\mathbb{R}^n}\left(M^{\sharp} f(x)\right)^p\omega(x)dx$$
for $1<p\le s_0$. Then, apply Lemma 3.5,
$$\|Mf\|_{(\mathbf{X}_1,\mathbf{X}_2)}\le C\|M^{\sharp} f\|_{(\mathbf{X}_1,\mathbf{X}_2)}$$
for $f\in\bigcup_{\omega\in A_1(\mathbb{R}^n)}L_{\omega}^s(\mathbb{R}^n)$. Hence, by Lemma 3.4
$$\|Mf\|_{(\mathbf{X}_1,\mathbf{X}_2)}\le C\|M^{\sharp} f\|_{(\mathbf{X}_1,\mathbf{X}_2)},~~f\in(\mathbf{X}_1,\mathbf{X}_2).$$
The proof is completed. $~~~~\blacksquare$

In fact, for the mixed-norm ball quasi-Banach function spaces $(\mathbf{X}_1,\mathbf{X}_2,\cdots,\mathbf{X}_m)$, the above discussions are right. Thus, we obtain the following corollary.

\textbf{Corollary 3.8} Let $\mathbf{X}_1$, $\mathbf{X}_2$,$\cdots$,$\mathbf{X}_m$ be a ball quasi-Banach function space and satisfy Assumption 3.2 for $1<s_0<\infty$. Then there exists a positive constant $C_0$ such that, for any $f\in(\mathbf{X}_1,\mathbf{X}_2,\cdots,\mathbf{X}_m)$,
$$\|Mf\|_{(\mathbf{X}_1,\mathbf{X}_2,\cdots,\mathbf{X}_m)}\le C\|M^{\sharp} f\|_{(\mathbf{X}_1,\mathbf{X}_2,\cdots,\mathbf{X}_m)}.\eqno{(3.6)}$$

Taking $m=1$, we get the following corollary.

\textbf{Corollary 3.9} Let $X$ be a ball quasi-Banach function space and satisfying Assumption 3.2 for $1<s_0<\infty$. Then there exists a positive constant $C_0$ such that, for any $f\in\mathbf{X}$,
$$\|Mf\|_{\mathbf{X}}\le C\|M^{\sharp} f\|_{\mathbf{X}}.\eqno{(3.7)}$$

If $1<\vec{p}=(p_1,p_2,\cdots,p_n)<\infty$, we define that $p_{-}=\min\{p_1,p_2,\cdots,p_n\}$ and $p_{+}=\max\{p_1,p_2,\cdots,p_n\}$.

By Remark 2.12, we know that mixed-norm Lebesgue spaces are a ball quasi-Banach function spaces. By Lemma 3.5 of \cite{24}, the $Mf$ is bounded on $L^{\vec{p}}(\mathbb{R}^n)$ with $1<\vec{p}=(p_1,p_2,\cdots,p_n)<\infty$. According to the dual theorem of  Theorem 1.a of \cite{7}, there exist $1<s<p_{-}$ such that $1<\frac{\vec{p}}{s}<\infty$ and the $Mf$ is bounded on $L^{(\frac{\vec{p}}{s})'}(\mathbb{R}^n)$. Hence, (3.7) holds when $\mathbf{X}=L^{\vec{p}}(\mathbb{R}^n)$. That is the following theorem.

\textbf{Theorem 3.10} Let $1<\vec{p}<\infty$. Then
$$\|Mf\|_{L^{\vec{p}}(\mathbb{R}^n)}\le C\|M^{\sharp} f\|_{L^{\vec{p}}(\mathbb{R}^n)}$$
hold for any $f\in{L^{\vec{p}}(\mathbb{R}^n)}$.

\textbf{Remark 3.11} (1) It is well-known that $M$ is bounded on classical Lebesgue spaces. Thus, by Corollary 3.8, Theorem 3.10 can be obtained more simply.

(2) Comparing Corollary 3.8 and Corollary 3.9, we find that if $Y$ is a mixed-norm ball Banach function space, then using Corollary 3.8 will be more simple. But Corollary 3.9 can be applied in a wider range.

\section{Application of Theorem 3.10}\label{sec3}
\par

In this section, a necessary and sufficient conditions of boundedness of commutator of $I_{\alpha}$ is given on mixed-norm Lebesgue spaces. Particularly, we point out that $\vec{p}\neq\vec{q}$ means there exist $i_0$ such that $p_{i_0}\neq q_{i_0}$.

In 2020, Zhang and Zhou gave necessary and sufficient conditions \cite{18}. Their result is stated as the following.

\textbf{Lemma 4.1} Let $0<\alpha<n,~1<\vec{p},\vec{q}<\infty$. Then
$$1<\vec{p}\le\vec{q}<\infty,~\vec{p}\neq\vec{q},~\alpha=\sum_{i=1}^n\frac{1}{p_i}-\sum_{i=1}^n\frac{1}{q_i}.$$
if and only if
$$\|I_{\alpha}f\|_{L^{\vec{q}}}\le C\|f\|_{L^{\vec{p}}}.$$

The following result can be proved by the means of Theorem 1.3 of \cite{23} and the fact $M_{\alpha}(f)(x)\lesssim I_{\alpha}(|f|)(x)$.

\textbf{Lemma 4.2} Let $0<\alpha<n$, $1<r<\infty$ and $b\in BMO(\mathbb{R}^n)$. Then there exists a constant $C>0$ independent of $b$ and $f$ such that
$$M^{\sharp}([b,I_{\alpha}](f))(x)\le C\|b\|_{BMO(\mathbb{R}^n)}\{I_{\alpha}(|f|)(x)+I_{r\alpha}(|f|^r)(x)^{\frac{1}{r}}\}.$$

The characterization of $BMO(\mathbb{R}^n)$ is given by the following theorem.

\textbf{Theorem 4.3} Let $0<\alpha<n,~1<\vec{p},\vec{q}<\infty$ and
$$1<\vec{p}\le\vec{q}<\infty,~\vec{p}\neq\vec{q},~\alpha=\sum_{i=1}^n\frac{1}{p_i}-\sum_{i=1}^n\frac{1}{q_i}.$$
Then, the following conditions are equivalent:\\
(a) $b\in BMO(\mathbb{R}^n)$.\\
(b) $[b,I_\alpha]$ is bounded from $L^{\vec{p}}$ to $L^{\vec{q}}$.

\textbf{Proof} (1)
By theorem 3.10, Lemma 4.2 and Lemma 4.1,
\begin{eqnarray*}
\|[b,I_{\alpha}]f\|_{L^{\vec{q}}(\mathbb{R}^n)}&\le&C\|M([b,I_{\alpha}](f))\|_{L^{\vec{q}}(\mathbb{R}^n)}\\
&\le&C\|M^{\sharp}([b,I_{\alpha}](f))\|_{L^{\vec{q}}(\mathbb{R}^n)}\\
&\le&C\|b\|_{BMO(\mathbb{R}^n)}(\|I_{\alpha}(|f|)\|_{L^{\vec{q}}(\mathbb{R}^n)} +\|I_{r\alpha}(|f|^r)^{\frac{1}{r}}\|_{L^{\vec{q}}(\mathbb{R}^n)})\\
&=&C\|b\|_{BMO(\mathbb{R}^n)}(\|I_{\alpha}(|f|)\|_{L^{\vec{q}}(\mathbb{R}^n)} +\|I_{r\alpha}(|f|^r)\|^{\frac{1}{r}}_{L^{\vec{q}/r}(\mathbb{R}^n)})\\
&\le&C\|b\|_{BMO(\mathbb{R}^n)}\|f\|_{L^{\vec{p}}(\mathbb{R}^n)}
\end{eqnarray*}

(2) Assume that $[b,I_\alpha]$ is bounded from $L^{\vec{p}}$ to $L^{\vec{q}}$. We use the same method as Janson \cite{25}. Choose $0\neq z_0\in\mathbb{R}^n$ such that $0\notin Q(z_0,2\sqrt{n})$. Then for $x\in Q(z_0,2\sqrt{n})$, $|x|^{n-\alpha}\in C^{\infty}(Q(z_0,2\sqrt{n}))$. Hence, $|x|^{n-\alpha}$ can be written as the absolutely convergent Fourier series:
$$|x|^{n-\alpha}\chi_{Q(z_0,2)}(x)=\sum_{m\in \mathbb{Z}^n}a_me^{2im\cdot x}\chi_{Q(z_0,2\sqrt{n})}(x)$$
with $\sum_{m\in \mathbb{Z}^n}|a_m|<\infty$.

For any $x_0\in\mathbb{R}^n$ and $t>0$, let $Q=Q(x_0,t)$ and $Q_{z_0}=Q(x_0+z_0t,t)$. Let $s(x)=\overline{sgn(\int_{Q'}(b(x)-b(y))dy)}$. Then
\begin{align*}
&\quad\frac{1}{|Q|}\int_Q|b(x)-b_{Q_{z_0}}|\\
&=\frac{1}{|Q|}\frac{1}{|Q_{z_0}|}\int_Q\left|\int_{Q_{z_0}}(b(x)-b(y))dy\right|dx\\
&=\frac{1}{|Q|}\frac{1}{|Q_{z_0}|}\int_Q\int_{Q_{z_0}}s(x)(b(x)-b(y))dydx.
\end{align*}
If $x\in Q$ and $y\in Q_{z_0}$, then $\frac{y-x}{t}\in Q(z_0,2\sqrt{n})$. Thereby,
\begin{align*}
&\quad\frac{1}{|Q|}\int_Q|b(x)-b_{Q_{z_0}}|\\
&=t^{-2n}\int_Q\int_{Q_{z_0}}s(x)(b(x)-b(y))|x-y|^{\alpha-n}|x-y|^{n-\alpha}dydx\\
&=t^{-2n}\int_Q\int_{Q_{z_0}}s(x)(b(x)-b(y))|x-y|^{\alpha-n}|x-y|^{n-\alpha}dydx\\
&=t^{-n-\alpha}\int_Q\int_{Q_{z_0}}s(x)(b(x)-b(y))|x-y|^{\alpha-n}\left(\frac{|x-y|}{t}\right)^{n-\alpha}dydx\\
&=t^{-n-\alpha}\sum_{m\in\mathbb{Z}^n}a_m\int_Q\int_{Q_{z_0}}s(x)(b(x)-b(y))|x-y|^{\alpha-n}e^{-2im\cdot \frac{y}{t}}dy\times e^{2im\cdot \frac{x}{t}}dx\\
&=t^{-n-\alpha}\sum_{m\in\mathbb{Z}^n}a_m\int_Q[b,I_\alpha](e^{-2im\cdot \frac{\cdot}{t}}\chi_{Q_{z_0}})(x)\times s(x)e^{2im\cdot \frac{x}{t}}dx.
\end{align*}
By H\"older for mixed-norm Lebesgue spaces,
$$\frac{1}{|Q|}\int_Q|b(x)-b_{Q_{z_0}}|
\le t^{-n-\alpha}\sum_{m\in\mathbb{Z}^n}a_m\|[b,I_\alpha](e^{-2im\cdot \frac{\cdot}{t}}\chi_{Q_{z_0}})\|_{L^{\vec{q}}(\mathbb{R}^n)} \|s\cdot e^{-2im\cdot \frac{\cdot}{t}}\chi_Q\|_{L^{\vec{q}'}(\mathbb{R}^n)}.$$
It is easy to calculate
$$\|s\cdot e^{-2im\cdot \frac{\cdot}{t}}\chi_Q\|_{L^{\vec{q}'}(\mathbb{R}^n)} =\|\chi_Q\|_{L^{\vec{q}'}(\mathbb{R}^n)}=t^{\sum_{i=1}^n\frac{1}{{q_i}'}}.$$
Hence,
$$\frac{1}{|Q|}\int_Q|b(x)-b_{Q_{z_0}}| =t^{-n-\alpha+\sum_{i=1}^n\frac{1}{{q_i}'}}\sum_{m\in\mathbb{Z}^n}a_m\|[b,I_\alpha](e^{-2im\cdot \frac{\cdot}{t}}\chi_{Q_{z_0}})\|_{L^{\vec{q}}(\mathbb{R}^n)}.$$
According to the hypothesis
\begin{align*}
&\quad\frac{1}{|Q|}\int_Q|b(x)-b_{Q_{z_0}}|\\
&\le t^{-n-\alpha+\sum_{i=1}^n\frac{1}{{q_i}'}}\sum_{m\in\mathbb{Z}^n}a_m\|e^{-2im\cdot \frac{\cdot}{t}}\chi_{Q_{z_0}}\|_{L^{\vec{p}}(\mathbb{R}^n)}\|[b,I_\alpha]\|_{L^{\vec{p}}(\mathbb{R}^n)\rightarrow L^{\vec{q}}(\mathbb{R}^n)}\\
&=t^{-n-\alpha+\sum_{i=1}^n\frac{1}{{q_i}'}+\sum_{i=1}^n\frac{1}{{p_i}}}\sum_{m\in\mathbb{Z}^n}a_m\|[b,I_\alpha]\|_{L^{\vec{p}}(\mathbb{R}^n)\rightarrow L^{\vec{q}}(\mathbb{R}^n)}\\
&\le \sum_{m\in\mathbb{Z}^n}|a_m|\|[b,I_\alpha]\|_{L^{\vec{p}}(\mathbb{R}^n)\rightarrow L^{\vec{q}}(\mathbb{R}^n)}\le C\|[b,I_\alpha]\|_{L^{\vec{p}}(\mathbb{R}^n)\rightarrow L^{\vec{q}}(\mathbb{R}^n)}.
\end{align*}
Thus, we have
$$\frac{1}{|Q|}\int_Q|b(x)-b(y)|dx\le\frac{2}{|Q|}\int_Q|b(x)-b_{Q_{z_0}}|dx\le C\|[b,I_\alpha]\|_{L^{\vec{p}}(\mathbb{R}^n)\rightarrow L^{\vec{q}}(\mathbb{R}^n)}$$
This prove $b\in BMO(\mathbb{R}^n)$. $~~~~\blacksquare$

\textbf{Remark 4.4} If $$\alpha=\sum_{i=1}^n\frac{1}{p_i}-\sum_{i=1}^n\frac{1}{q_i},~p_j\sum_{i=1}^n\frac{1}{p_i}=q_j\sum_{i=1}^n\frac{1}{q_i}~~(j=1,\cdots,n),$$
it easy to prove
$$\sum_{i=1}^n\frac{1}{p_i}>\sum_{i=1}^n\frac{1}{q_i},~1<\vec{p}<\vec{q}<\infty.$$

\section{Characterization of homogeneous Lipschitz space}\label{sec3}
\par
In this section, a characterization of homogeneous Lipschitz spaces is given. Let us recall the definition of homogeneous Lipschitz spaces.

\textbf{Definition 5.1} Let $0<\beta<1$. The definition of homogeneous Lipschitz space is defined by
$$\dot{\Lambda}_{\beta}:=\{f:|f(x)-f(y)|\le C|x-y|^{\beta}\}.$$

The following lemma can be found in Lemma 1.5 of \cite{26}.

\textbf{Lemma 5.2} If $0<\beta<1$ and $1<q\le\infty$, then
\begin{align*}
\|f\|_{\dot{\Lambda}_{\beta}}&\approx\sup_{Q\subset\mathbb{R}^n}\frac{1}{|Q|^{1+\beta/n}}\int_{Q}|f(y)-f_Q|dy\\
&\approx\sup_{Q\subset\mathbb{R}^n}\frac{1}{|Q|^{\beta/n}}\left(\frac{1}{|Q|}\int_{Q}|f(y)-f_Q|^qdy\right)^{\frac{1}{q}},
\end{align*}
for $q=\infty$ the formula should be interpreted appropriately, where the supremum is taken over all cubes $Q$ in $\mathbb{R}^n$.

\textbf{Theorem 5.3} Let $0<\alpha<n,~0<\beta<1,~1<\vec{p},\vec{q}<\infty$ and
$$1<\vec{p}\le\vec{q}<\infty,~\vec{p}\neq\vec{q},~\alpha+\beta=\sum_{i=1}^n\frac{1}{p_i}-\sum_{i=1}^n\frac{1}{q_i}.$$
Then, the following conditions are equivalent:\\
(a) $b\in \dot{\Lambda}_{\beta}$.\\
(b) $[b,I_\alpha]$ is bounded from $L^{\vec{p}}(\mathbb{R}^n)$ to $L^{\vec{q}}(\mathbb{R}^n)$.

\textbf{Proof} (1) Let $b\in \dot{\Lambda}_{\beta}$. Then
\begin{align*}
|[b,I_{\alpha}]f(x)|&=\left|\int_{\mathbb{R}^n}\frac{(b(x)-b(y))f(y)}{|x-y|^{n-\alpha}}dy\right|\\
&\le\int_{\mathbb{R}^n}\frac{|b(x)-b(y)|\cdot |f(y)|}{|x-y|^{n-\alpha}}dy\\
&\le\|b\|_{\dot{\Lambda}_{\beta}}\int_{\mathbb{R}^n}\frac{|f(y)|}{|x-y|^{n-(\alpha+\beta)}}dy\\
&=\|b\|_{\dot{\Lambda}_{\beta}}I_{\alpha+\beta}(|f|)(x).
\end{align*}

According to Lemma 4.1,
$$\|[b,I_{\alpha}]f\|_{L^{\vec{q}}(\mathbb{R}^n)}\le C\|b\|_{\dot{\Lambda}_{\beta}}\|I_{\alpha+\beta}(|f|)\|_{L^{\vec{q}}(\mathbb{R}^n)}\le C\|b\|_{\dot{\Lambda}_{\beta}}\|f\|_{L^{\vec{p}}(\mathbb{R}^n)}.$$

(2) Let $[b,I_\alpha]$ is bounded from $L^{\vec{p}}(\mathbb{R}^n)$ to $L^{\vec{q}}(\mathbb{R}^n)$. Let $Q$, $Q_{z_0}$ and $s(x)$ is the same as Theorem 4.4. It is easy to calculate that
\begin{align*}
&\quad\frac{1}{|Q|}\int_Q|b(x)-b_{Q_{z_0}}|\\
&=t^{-n-\alpha}\sum_{m\in\mathbb{Z}^n}a_m\int_Q[b,I_\alpha](e^{-2im\cdot \frac{\cdot}{t}}\chi_{Q_{z_0}})(x)\times s(x)e^{2im\cdot \frac{x}{t}}dx\\
&\le t^{-n-\alpha}\sum_{m\in\mathbb{Z}^n}a_m\|[b,I_\alpha](e^{-2im\cdot \frac{\cdot}{t}}\chi_{Q_{z_0}})\|_{L^{\vec{q}}(\mathbb{R}^n)} \|s\cdot e^{-2im\cdot \frac{\cdot}{t}}\chi_Q\|_{L^{\vec{q}'}(\mathbb{R}^n)}\\
&=t^{-n-\alpha+\sum_{i=1}^n\frac{1}{{q_i}'}+\sum_{i=1}^n\frac{1}{{p_i}}}\sum_{m\in\mathbb{Z}^n}a_m\|[b,I_\alpha]\|_{L^{\vec{p}}(\mathbb{R}^n)\rightarrow L^{\vec{q}}(\mathbb{R}^n)}\\
&\le t^{\beta}\sum_{m\in\mathbb{Z}^n}|a_m|\|[b,I_\alpha]\|_{L^{\vec{p}}(\mathbb{R}^n)\rightarrow L^{\vec{q}}(\mathbb{R}^n)}\le Ct^{\beta}\|[b,I_\alpha]\|_{L^{\vec{p}}(\mathbb{R}^n)\rightarrow L^{\vec{q}}(\mathbb{R}^n)}.
\end{align*}
Therefore,
$$\frac{1}{|Q|^{1+\beta/n}}\int_Q|b(x)-b_{Q}|\le \frac{2}{|Q|^{1+\beta/n}}\int_Q|b(x)-b_{Q_{z_0}}|\le 2C\|[b,I_{\alpha}]\|_{L^{\vec{p}}(\mathbb{R}^n)\rightarrow L^{\vec{q}}(\mathbb{R}^n)}.$$
Due to Lemma 5.2, the proof is completed. $~~~~\blacksquare$

By the process of Theorem 5.3, the following corollary holds.

\textbf{Corollary 5.4} Let $Tf$ is a operator and its commutator $[b,T](f)$ satisfies that
$$|[b,T]f(x)|\le C  I_{\alpha,b}(|f|)(x).$$
Let $0<\alpha<n,~0<\beta<1,~1<\vec{p},\vec{q}<\infty$. If
$$1<\vec{p}\le\vec{q}<\infty,~\vec{p}\neq\vec{q},~\alpha+\beta=\sum_{i=1}^n\frac{1}{p_i}-\sum_{i=1}^n\frac{1}{q_i},$$
and $b\in \dot{\Lambda}_{\beta}$, then
$$\|[b,T]f\|_{L^{\vec{q}}(\mathbb{R}^n)}\le C\|b\|_{\dot{\Lambda}_{\beta}}\|f\|_{L^{\vec{p}}(\mathbb{R}^n)}.$$

\textbf{Proof} Let $b\in \dot{\Lambda}_{\beta}$. Then
\begin{align*}
|[b,T]f(x)|&\le C\int_{\mathbb{R}^n}\frac{|b(x)-b(y)|\cdot |f(y)|}{|x-y|^{n-\alpha}}dy\\
&\le C\|b\|_{\dot{\Lambda}_{\beta}}\int_{\mathbb{R}^n}\frac{|f(y)|}{|x-y|^{n-(\alpha+\beta)}}dy\\
&= C\|b\|_{\dot{\Lambda}_{\beta}}I_{\alpha+\beta}(|f|)(x).
\end{align*}
According to Lemma 4.1,
$$\|[b,I_{\alpha}]f\|_{L^{\vec{q}}(\mathbb{R}^n)}\le C\|b\|_{\dot{\Lambda}_{\beta}}\|I_{\alpha+\beta}(|f|)\|_{L^{\vec{q}}(\mathbb{R}^n)}\le C\|b\|_{\dot{\Lambda}_{\beta}}\|f\|_{L^{\vec{p}}(\mathbb{R}^n)}.$$
The proof is completed. $~~~~\blacksquare$

The corollary is very useful and two examples are given in the following section.

\section{Two applications of Corollary 5.4}\label{sec3}
\par

\textbf{Example 6.1} The fractional maximal function is defined as
$$M_{\alpha}f(x):=\sup_{Q\ni x}\frac{1}{|Q|^{1-\frac{\alpha}{n}}}\int_{Q}|f(y)|dy,$$
where the supremum is taken over all cube $Q\subset\mathbb{R}^n$ containing $x$ and its commutator is defined by
$$M_{\alpha,b}f(x):=\sup_{Q\ni x}\frac{1}{|Q|^{1-\frac{\alpha}{n}}}\int_{Q}|b(x)-b(y)||f(y)|dy,$$
where $b$ a locally integrable function. It is easy to prove that
$$|M_{\alpha,b}f(x)|\le C I_{\alpha,b}(|f|)(x).$$

Before the second example, let us recall generalized fractional integral operators.

Suppose that $\mathcal{L}$ is a linear operator which generates an analytic semigroup $\{e^{-t\mathcal{L}}\}_{t>0}$ on
$L^2(\mathbb{R}^n)$ with a kernel $p_t(x,y)$ satisfying Gaussian upper bound; that is
$$|p_t(x,y)|\le\frac{C_1}{t^{n/2}}e^{-C_2\frac{|x-y|^2}{t}}~~x,y\in \mathbb{R}^n,$$
where $C_1,C_2>0$ are independent of $x,~y$ and $t$.

For any $0<\alpha<n$, the generalized fractional integrals $\mathcal{L}^{-\alpha/2}$ associated with the operator $\mathcal{L}$ is defined by
$$\mathcal{L}^{-\alpha/2}f(x)=\frac{1}{\Gamma(\alpha/2)}\int_0^\infty e^{-t\mathcal{L}}(f)(x)\frac{dt}{t^{-\alpha/2+1}}.$$

Note that if $\mathcal{L}=-\Delta$ is the Laplacian on $\mathbb{R}^n$, then $\mathcal{L}^{-\alpha/2}$ is the classical fractional intergral operator $I_{\alpha}$. See, for example, Chapter 5 of \cite{27}. Since the semigroup $\mathcal{L}^{-\alpha/2}$ has a kernel $p_t(x,y)$, it is easy to check that for all $x\in\mathbb{R}^n$£¬
$$|\mathcal{L}^{-\alpha/2}f(x)|\le CI_{\alpha}(|f|)(x).$$
(see \cite{28}). In fact, if we denote the the kernel of $\mathcal{L}^{-\alpha/2}$ by $K_{\alpha}(x,y)$, it is easy to obtain that
\begin{align*}
\mathcal{L}^{-\alpha/2}f(x)&=\frac{1}{\Gamma(\alpha/2)}\int_0^\infty e^{-t\mathcal{L}}(f)(x)\frac{dt}{t^{-\alpha/2+1}}\\
&=\frac{1}{\Gamma(\alpha/2)}\int_0^\infty\int_{\mathbb{R}^n}p_t(x,y)f(y)dy\frac{dt}{t^{-\alpha/2+1}}\\
&=\int_{\mathbb{R}^n}\frac{1}{\Gamma(\alpha/2)}\int_0^\infty p_t(x,y)\frac{dt}{t^{-\alpha/2+1}}\cdot f(y)dy\\
&=\int_{\mathbb{R}^n}K_{\alpha}(x,y)\cdot f(y)dy.
\end{align*}
Hence, by Gaussian upper bound,
\begin{align*}
|K_{\alpha}(x,y)|&=\left|\frac{1}{\Gamma(\alpha/2)}\int_0^\infty p_t(x,y)\frac{dt}{t^{-\alpha/2+1}}\right|\\
&\le \frac{1}{\Gamma(\alpha/2)}\int_0^\infty |p_t(x,y)|\frac{dt}{t^{-\alpha/2+1}}\\
&\le C\int_0^\infty e^{-C_2\frac{|x-y|^2}{t}}\frac{dt}{t^{n/2-\alpha/2+1}}\\
&\le C\cdot\frac{1}{|x-y|^{n-\alpha}}
\end{align*}

\textbf{Example 6.2} Let $b$ a locally integrable function. If the commutators of generalized fractional integral operators generated by $b$ and $\mathcal{L}^{-\alpha/2}$ are defined by
$$[b,\mathcal{L}^{-\alpha/2}]:=b(x)\mathcal{L}^{-\alpha/2}(f)(x)-\mathcal{L}^{-\alpha/2}(bf)(x).$$
It is obvious that
\begin{align*}
|[b,\mathcal{L}^{-\alpha/2}]f(x)|&=\left|\int_{\mathbb{R}^n}(b(x)-b(y))K_{\alpha}(x,y)f(y)dy\right|\\
&\le\int_{\mathbb{R}^n}|b(x)-b(y)||K_{\alpha}(x,y)||f(y)|dy\\
&\le C\int_{\mathbb{R}^n}|b(x)-b(y)|\frac{|f(y)|}{|x-y|^{n-\alpha}}dy\\
&=CI_{\alpha,b}(|f|)(x).
\end{align*}


\begin{thebibliography}{1}
\bibitem{1} Stein E. M., Singular Integrals and Differentiability of Functions. Princeton University Press, Princeton (1970)
\bibitem{2} Stein E. M. and Weiss G., Introduction to Fourier analysis on Euclidean spaces. Princeton University Press, Princeton (1971)
\bibitem{3} Antonic N. and Ivec I., On the H\"{o}rmander-Mihlin theorem for mixed-norm Lebesgue spaces, Math. Anal. Appl. 433:176-199 (2016)
\bibitem{4}  Kenig C. E., On the local and global well-posedness theory for the KP-I equation, Ann. Inst. H. Poincar Anal. Non Linaire 21:827-838 (2004)
\bibitem{5} Kim D., Elliptic and parabolic equations with measurable coefficients in $L^{p}$-spaces with mixed norms, Methods Appl. Anal. 15:437-468 (2008)
\bibitem{6} Krylov N.V., Parabolic equations with VMO coefficients in Sobolev spaces with mixed norms, Funct. Anal. 250:521-558 (2007)
\bibitem{7} Benedek A. and Panzone R., The space $L^{p}$, with mixed norm, Duke Math. 28:301-324 (1961)
\bibitem{8} Fernandez D.L., Lorentz spaces, with mixed norms, Funct. Anal. 25:128-146 (1977)
\bibitem{9} Milman M., Embeddings of Lorentz-Marcinkiewicz spaces with mixed norms, Anal. Math. 4:215-223 (1978)
\bibitem{10} Milman M., A note on L(p,q) spaces and Orlicz spaces with mixed norms, Proc. Amer. Math. Soc. 83:743-746 (1981)
\bibitem{11} Cleanthous G., Georgiadis A.G. and Nielsen M., Anisotropic mixed-norm Hardy spaces, J. Geom. Anal. (2017)
\bibitem{12} Besov O. V., Il'in V. P. and Nikolski{\i} S.M., Integral Representations of Functions, and Embedding Theorems, Second edition, Fizmatlit "Nauka", Moscow, 480,  (Russian,1996)
\bibitem{13} Chen T. and  Sun W., Iterated and Mixed Weak Norms with Applications to Geometric Inequalities. (2017) DOI: 10.1007/s12220-019-00243-x
\bibitem{14} Huang L. and Yang D., On Function Spaces with Mixed Norms -- A Survey. (2019) arXiv:1908.03291v1
\bibitem{15} Nogayama T., Mixed Morrey spaces, Positivity (2019) DOI: 10.1007/s11117-019-00646-8
\bibitem{16} Nogayama T., Boundedness of commutators of fractional integral operators on mixed Morrey spaces, Integral Transform. Spec. Funct. (2019)
\bibitem{17} Chanillo S., A note on commutators. Indiana Univ. Math. 31(1):7-16 (1982)
\bibitem{18} Zhang H. and Zhou J.£¬ The Boundedness of Fractional Integral Operators in Local and Global Mixed Morrey-type Spaces. (2021) arXiv:2102.01304v1
\bibitem{19} Sawano Y., Ho K.-P., Yang D. and Yang S., Hardy spaces for ball quasi-Banach function spaces, Diss. Math. 525:1-102 (2017)

\bibitem{21} Chang D.C., Wang S., Yang D. and Zhang Y., Littlewood-Paley characterizations of Hardy-type spaces associated with ball quasi-Banach function spaces, Complex Anal. Oper. Theory 14 (2020), Paper No. 40, 33 pp.
\bibitem{22} Bennett C. and Sharpley R., Interpolation of Operators, Pure Appl. Math. 129, Academic Press, Boston, MA, 1988.
\bibitem{23}  Cruz-Uribe D. and Fiorenza A., Endpoint estimates and weighted norm inequalities for commutators of fractional integrals. Publ. Mat. 47:103-131 (2003)
\bibitem{24} Huang L., Liu J., Yang D. and Yuan W., Atomic and Littlewood-Paley characterizations of anisotropic mixed-norm Hardy spaces and their applications, J. Geom. Anal. 29:1991-2067 (2019)
\bibitem{25} Janson S., Mean oscillation and commutators of singular integral operators. Ark. Math. 16:263-270 (1978)
\bibitem{26} Paluszy\'nski M., Characterization of the Besov spaces via the commutator operator of Coifman, Rochberg and Weiss. Indiana Univ. Math. J. 44:1-17 (1995)
\bibitem{27} Stein, E.M., Singular Integrals and Differentiability of Functions. Princeton University Press, Princeton (1970)
\bibitem{28} Duong, X.T., Yan, L.X., On commutators of fractional integrals. Proc. Am. Math. Soc. 132(12):3549-3557 (2004)
\bibitem{29} Garc¨ªa-Cuerva and Francia J., Weighted Norm Inequalities and Related Topics. Elsevier, London (1985)
\bibitem{30} Ho K.-P., Strong maximal operator on mixed-norm spaces. Ann. Univ. Ferrara, 62(2):1-17 (2016)
\bibitem{31} Morrey C. B., On the solutions of quasi-linear elliptic partial differential equations, Trans. Amer. Math. Soc. 43:126-166 (1938)
\bibitem{32} O. Kov\'a\u{a}ik and J. R\'akosn\'ik, On spaces $L^{p(x)}$ and $W^{k,p(x)}$, Czechoslovak Math. 41(116):592-618 (1991)
\bibitem{33} Birnbaum Z. and Orlicz W., \"Uber die verallgemeinerung des begriffes der zueinander konjugierten potenzen, Stud. Math. 3:1-67 (1931)
\bibitem{34} Grafakos L.. Modern Fourier Analysis. Springer New York, (2009)
\bibitem{35} Zhang Y., Yang D., Yuan W. and Wang S., Real-variable characterizations of Orlicz-slice Hardy spaces, Anal. Appl. (Singap.) 17:597-664 (2019)
\bibitem{36} Bennett C. and Sharpley R.£¬ Interpolation of Operators. Academic Press, London (1988)
\bibitem{37} Sawano Y., Theory of Besov Spaces, Developments in Mathematics 56, Springer, Singapore, (2018)
\bibitem{38} Izuki M. and Sawano Y., Characterization of BMO via ball Banach function spaces, Vestn. St.-Peterbg. Univ. Mat. Mekh. Astron. 4(62):78-86 (2017)
\bibitem{39} Izuki M., Noi T. and Sawano Y., The John-Nirenberg inequality in ball Banach function spaces and application to characterization of BMO, J. Inequal. Appl. 2019(1) (2019)
\bibitem{40} Wang S., Yang D., Yuan W. and Zhang Y., Weak Hardy-type spaces associated with ball quasi-Banach function spaces II: Littlewood-Paley characterizations and real interpolation, J. Geom. Anal. (2019), DOI: 10.1007/s12220-019-00293-1.
\bibitem{41} Zhang Y., Yang D., Yuan W. and Wang S., Real-variable characterizations of Orlicz-slice Hardy spaces, Anal. Appl. (Singap.) 17:597-664 (2019)
\bibitem{42} Chanillo S., A note on commutators, Indiana Univ. Math. J. 31:7-16 (1982)
\bibitem{43} Paluszynski M., Characterization of the Besov spaces via the commutator operator of Coifman, Rochberg and Weiss, Indiana Univ. Math. J. 44:1-17 (1995)
\bibitem{44} Blozinski A., Multivariate rearrangements and Banach function spaces with mixed norms. Trans. Am. Math. Soc. 263: 149-167 (1981)
\end{thebibliography}
\end{document}